\documentclass[a4paper,11pt,twoside]{amsart}
\usepackage[english]{babel} 
\usepackage[utf8]{inputenc} 
\usepackage[T1]{fontenc} 

\usepackage{amssymb,amsmath,amsfonts,amsxtra,amsthm} 
\usepackage{stmaryrd} 
\usepackage{mathtools} 
\usepackage{mathrsfs} 
\usepackage{cases} 
\usepackage[all]{xy} 

\usepackage{geometry} 
\usepackage{lmodern} 
\usepackage{xargs} 
\usepackage{graphicx} 
\usepackage[svgnames,table]{xcolor} 
\usepackage[colorlinks=true,citecolor=Magenta,linkcolor=Blue, urlcolor=Green, pdftoolbar=true]{hyperref}
\usepackage{enumitem}
\usepackage{etaremune} 

\newcommand{\ensemblenombre}[1]{\ensuremath{\mathbb{#1}}}
\newcommand{\N}{\ensemblenombre{N}}
\newcommand{\Z}{\ensemblenombre{Z}}

\newcommand{\R}{\ensemblenombre{R}}

\newcommand{\intervalle}[4]{\ensuremath{\mathopen{#1}#2
		\mathclose{}\mathpunct{};#3
		\mathclose{#4}}}
\newcommand{\intervalleff}[2]{\intervalle{[}{#1}{#2}{]}}

\newcommand{\piun}[1]{\pi_1({#1})}

\newcommand{\RP}[1]{\ensuremath{\R\mathbf{P}^{#1}}}
\newcommand{\Sn}[1]{\mathbf{S}^{#1}}

\newcommand{\Tn}[1]{\mathbf{T}^{#1}}

\newcommand{\Falpha}{\ensuremath{\mathcal{F}_\alpha}}
\newcommand{\Falphazero}{\ensuremath{\mathcal{F}_\alpha^0}}
\newcommand{\Fbeta}{\ensuremath{\mathcal{F}_\beta}}
\newcommand{\Fa}{\ensuremath{F_\alpha}}
\newcommand{\Fz}{\ensuremath{F_\alpha^0}}

\newcommand{\PSL}[1]{\ensuremath{\mathrm{PSL}_{#1}\mathopen{(}\R\mathclose{)}}}

\newcommand{\Deftheta}{\ensuremath{\mathsf{Def}_{\theta}(\Tn{2},\mathsf{0})}}

\newcommand{\Homologie}[3]{\ensuremath{\mathrm{H}_{#1}\mathopen{(}#2\mathpunct,#3\mathclose{)}}}

\DeclareMathOperator{\id}{id}

\DeclareMathOperator{\Homeo}{Homeo}

\numberwithin{equation}{section}
\numberwithin{figure}{section}

\AtBeginDocument{}

\addto\captionsenglish{}

\theoremstyle{definition}
\newtheorem{definition}{Definition}[section]
\newtheorem{propositiondefinition}[definition]{Proposition-Definition}

\theoremstyle{plain}

\newtheorem{theoremintro}{Theorem}

\newtheorem{corollaryintro}[theoremintro]{Corollary}

\newtheorem{lemma}[definition]{Lemma}
\newtheorem{proposition}[definition]{Proposition}

\theoremstyle{remark}

\title{On simultaneous conjugacies of
pairs of transverse foliations of the torus}
\author{Martin Mion-Mouton}
\date{\today}

\begin{document}
\address{
Martin Mion-Mouton,
Max Planck Institute for Mathematics in the Sciences
in Leipzig
}
\email{\href{mailto:martin.mion@mis.mpg.de}{martin.mion@mis.mpg.de}}
\urladdr{\url{https://orcid.org/0000-0002-8814-0918}}

\begin{abstract}
We prove in this note that
two pairs of transverse minimal topological foliations of the torus
are individually conjugated if, and only if
they are simultaneously conjugated.
\end{abstract}

\maketitle

\section{Introduction}
It is known since Poincaré
\cite{poincare_sur_1885}
that the \emph{rotation number}
(introduced in Proposition-Definition \ref{propositiondefinition-nombrerotation} below) is the only invariant of minimal,
orientation-preserving
homeomorphisms of the circle $\Sn{1}=\R/\Z$
under topological conjugacies.
For a one-dimensional oriented topological foliation
$\mathcal{F}$ of the torus
$\Tn{2}=\R^2/\Z^2$,
the global conjugacy
invariant which naturally replaces the rotation number is
the \emph{projective asymptotic cycle}
due to Schwartzman \cite{schwartzman_asymptotic_1957}
(see \eqref{equation-asymptoticcycle}).
The latter is a half-line
$A(\mathcal{F})$ in the first homology group
$\Homologie{1}{\Tn{2}}{\R}$ of the torus,
which gives a rigorous idea to
the ``asymptotic direction of $\mathcal{F}$ in homology'',
and reflects the dynamics of the foliation in multiple ways.
For instance, the asymptotic cycle of a minimal
foliation of $\Tn{2}$ is \emph{irrational}
in the sense that the half-line
$A(\mathcal{F})$
does not intersect
$\piun{\Tn{2}}\equiv\Z^2\subset\R^2\equiv\Homologie{1}{\Tn{2}}{\R}$.
Moreover in the same way than the rotation number for minimal circle homeomorphisms,
the projective asymptotic cycle is known to be a complete
conjugacy invariant for minimal foliations.
Any minimal oriented topological foliation of $\Tn{2}$
is indeed conjugated
to the linear foliation defined by its asymptotic cycle
$A(\mathcal{F})$, by a homeomorphism isotopic to the identity
(see Proposition \ref{proposition-classificationfeuilletagesminimaux}).
\par We study in this note the \emph{minimal topological bi-foliations}
of the torus $\Tn{2}=\R^2/\Z^2$,
defined as the pairs $(\Falpha,\Fbeta)$ of minimal oriented
topological foliations
of $\Tn{2}$
which are
\emph{transverse} (in a sense made precise
in Definition \ref{definition-transverse}).
Such objects appear naturally in different dynamical contexts.
For instance, any closed surface $S$ which is a transverse section
of a three-dimensional Anosov flow is endowed
with a bi-foliation, whose leaves are the intersections
of $S$ with the weak stable and unstable leaves
of the Anosov flow.
While suspensions are the only three-dimensional Anosov flows
admitting such a transverse section,
Birkhoff sections \cite{birkhoffDynamicalSystemsTwo1917}
constitute a natural generalization of them
which is an important tool for
the classification of Anosov flows,
and any genus one Birkhoff section of a transitive
three-dimensional Anosov flow inherits in the same way
a minimal topological bi-foliation.
In a recent work of the author
\cite{mion-moutonRigiditySingularDeSitter2024},
bi-foliations of the torus are studied from a geometrical point of view
in the setting of
\emph{singular de-Sitter Lorentzian metrics}.
Isotopy classes of singular de-Sitter metrics of $\Tn{2}$
are indeed shown to be equivalent to the
\emph{simultaneous} conjugacy classes
of their minimal lightlike bi-foliations
(see \cite[Theorem A]{mion-moutonRigiditySingularDeSitter2024}).
Two bi-foliations
$(\Falpha^1,\Fbeta^1)$ and $(\Falpha^2,\Fbeta^2)$
are said to be \emph{simultaneously
conjugated} by a homeomorphism
$\phi$ of $\Tn{2}$,
if $\Falpha^2=\phi_*\Falpha^1$ and $\Fbeta^2=\phi_*\Fbeta^1$.
\par If two bi-foliations $(\Falpha^1,\Fbeta^1)$
and $(\Falpha^2,\Fbeta^2)$ of the torus
have the same pairs of projective asymptotic cycles, then
we already know
that $\Falpha^1$ and $\Falpha^2$ in the one hand,
and that $\Fbeta^1$ and $\Fbeta^2$ in the other hand,
are individually conjugated.
It is then natural to ask
if $(\Falpha^1,\Fbeta^1)$ and $(\Falpha^2,\Fbeta^2)$
are actually \emph{simultaneously} conjugated.
The main goal of the present note is
to answer this question positively
with the following result.
\begin{theoremintro}\label{theoremintro-bifoliations}
Let $(\Falpha,\Fbeta)$ be a minimal topological
bi-foliation of the torus
and $p\in\Tn{2}$.
Then $(\Falpha,\Fbeta)$ is simultaneously conjugated
to the linear bi-foliation
$(\mathcal{F}_{A(\Falpha)},\mathcal{F}_{A(\Fbeta)})$
defined by its asymptotic cycles,
by a homeomorphism isotopic to the identity relatively to $p$.
\end{theoremintro}
I was informed, after posting this preprint on arXiv,
that Theorem \ref{theoremintro-bifoliations}
was already proved in \cite{aransonClassificationSupertransitive2Webs2003}.
\par Theorem \ref{theoremintro-bifoliations} describes thus
a ``\emph{dynamical} Teichmüller space'':
the space of bi-foliations of the torus modulo simultaneous conjugacies
isotopic to the identity.
More precisely, it identifies the subset
of minimal bi-foliations in this Teichmüller space
(which are the most dynamically relevant ones)
with the subset $\{l_\alpha\neq l_\beta\}$ of $(\RP{1})^2$,
which is particularly useful
to study the dynamics of the
\emph{mapping class group} $\mathrm{Mod}^+(\Tn{2})$ of the torus
on its space of bi-foliations.
We recall that $\mathrm{Mod}^+(\Tn{2})$ is the quotient of
the group of homeomorphisms of $\Tn{2}$
by the subgroup of homeomorphisms isotopic to the identity,
which acts naturally on the space of bi-foliations
of $\Tn{2}$ modulo isotopies.
Theorem \ref{theoremintro-bifoliations} intertwines thus
(in restriction to minimal bi-foliations and irrational lines)
the \emph{a priori} complicated dynamics of $\mathrm{Mod}^+(\Tn{2})$
on the space of bi-foliations
with
the \emph{explicit} diagonal action of
$\mathrm{PSL}_2(\Z)$ on $(\RP{1})^2$.
The latter dynamics being entirely known,
one can hope to use this description
to deduce new informations.
\par This idea applies for instance to
the \emph{geometrical} Teichmüller space
$\Deftheta$ of
singular de-Sitter metrics of $\Tn{2}$
studied in
\cite{mion-moutonRigiditySingularDeSitter2024}.
It allows to reformulate its main result to
identify a part of $\Deftheta$
with the subset $\{l_\alpha\neq l_\beta\text{~irrational}\}$
of $(\RP{1})^2$, which intertwines
(on relevant subsets)
the action of $\mathrm{Mod}^+(\Tn{2})$ on $\Deftheta$
with the diagonal action of $\mathrm{PSL}_2(\Z)$
on $(\RP{1})^2$.
In a joint work in progress, Florestan Martin-Baillon and the author
use this description to study the dynamics
of $\mathrm{Mod}^+(\Tn{2})$ on relative
character varieties of
the one-holed torus with values in $\PSL{2}$.

As a direct consequence of Theorem \ref{theoremintro-bifoliations},
we obtain the following result by applying
an argument due to Ghys-Sergiescu
\cite{ghysStabiliteConjugaisonDifferentiable1980}.
\begin{corollaryintro}\label{corollaryintro-automorphismesbifoliations}
Let $f$ be a homeomorphism of $\Tn{2}$
 which is isotopic to the identity relatively to a point,
 and preserves a minimal topological bi-foliation.
 Then $f$ is the identity.
\end{corollaryintro}
Corollary \ref{corollaryintro-automorphismesbifoliations}
strengtens
the main result of \cite{mion-moutonRigiditySingularDeSitter2024},
by showing that any \emph{topological} conjugacy between the
minimal lightlike bi-foliations of two singular de-Sitter metrics
(having a unique singularity of the same angle)
is actually an \emph{isometry}.
This implies in particular a geometric rigidity
result for this class of dynamical systems:
any topological conjugacy between
such minimal lightlike bi-foliations
is actually piecewise smooth.

\section{Preliminaries}\label{section-definitions}

\subsection{Foliations of the torus and suspensions}
\label{subsection-foliationsasymptoticcycle}
We recall the basic definitions of foliations,
refering to
\cite[Chapter \MakeUppercase{\romannumeral 2}]{camacho_geometric_1985}
and \cite[Chapter \MakeUppercase{\romannumeral 1}]{hector_introduction_1986}
for more details.
\begin{definition}\label{definition-foliation}
Let $S$ be a topological surface.
A \emph{foliated (topological) atlas} of $S$
is a continuous atlas
$\mathcal{A}$ of $S$ satisfying the following conditions.
\begin{enumerate}
 \item For any chart $(U,\varphi)\in\mathcal{A}$,
 $\varphi(U)=I\times J\subset\R^2$
 with $I$ and $J$ open intervals in $\R$.
 \item For any $(U_i,\varphi_i),(U_j,\varphi_j)\in\mathcal{A}$
 such that $U_i\cap U_j\neq\varnothing$,
 the transition map
$\varphi_{i,j}\coloneqq\varphi_i\circ\varphi_j^{-1}\colon
\varphi_j(U_i\cap U_j)\to\varphi_i(U_i\cap U_j)$ is of the form
\begin{equation*}
 \varphi_{i,j}(x,y)=(\alpha_{i,j}(x,y),\gamma_{i,j}(y)).
\end{equation*}
\end{enumerate}
A \emph{(topological) foliation} $\mathcal{F}$
of $S$ is a maximal
foliated topological atlas of $S$,
the charts of which are called the \emph{foliation charts}
of $\mathcal{F}$.
In a domain $U$ of foliated chart of $\mathcal{F}$
and for any $p=\varphi^{-1}(x,y)\in U$,
$P^U_{\mathcal{F}}(p)=\varphi^{-1}(I\times\{y\})$
is called the \emph{plaque} of $x$ for $\mathcal{F}$.
The set of all plaques is a basis of a topology on $S$
whose connected components are called the \emph{leaves}
of $\mathcal{F}$, and
the leaf containing $x$ is denoted by $\mathcal{F}(x)$.
$\mathcal{F}$ is said \emph{oriented} by the choice of
a sub-atlas of $\mathcal{A}$ such that
$\alpha_{i,j}(\cdot,y)$ is an orientation-preserving map
for any $(i,j)$ and $y$.
In this case, each leaf of $\mathcal{F}$ inherits the orientation
given in any foliated chart by the identification
of the plaque $\varphi^{-1}(I\times\{y\})$
with $I\subset\R$.
Henceforth, all the foliations will implicitly be topological and oriented.
\par A \emph{topological conjugacy}
between two topological foliations
$(S_1,\mathcal{F}_1)$ and $(S_2,\mathcal{F}_2)$
is a homeomorphism $\phi\colon S_1\to S_2$
such that $\mathcal{F}_2=\phi_*\mathcal{F}_1$, namely
$\phi(\mathcal{F}_1(x))=\mathcal{F}_2(x)$
for any $x\in S_1$.
\end{definition}

Let $\mathcal{F}$ be a foliation of $\Tn{2}$
which admits a closed \emph{section},
\emph{i.e.} a simple closed curve
$\gamma\subset\Tn{2}$
transverse to $\mathcal{F}$
(in the sense of Definition \ref{definition-transverse} below)
and intersecting all of its leaves.
Then the \emph{first-return map}
\begin{equation*}\label{equation-firstreturn}
 P^\gamma_{\mathcal{F}}\colon\gamma\to\gamma
\end{equation*}
of $\mathcal{F}$ on $\gamma$
is well-defined,
$P^\gamma_{\mathcal{F}}(x)$ being the first intersection point
of the oriented leaf $\mathcal{F}(x)$ with $\gamma$ after $x$.
We can now describe $\mathcal{F}$ in terms of $P^\gamma_{\mathcal{F}}$
in the following way.
Given an orientation-preserving
homeomorphism $T$ of the circle,
the \emph{suspension}
$\mathcal{F}_T$ of $T$ is
the oriented foliation
of
the topological torus
\begin{equation*}
M_T\coloneqq \intervalleff{0}{1}\times\Sn{1}/\{(1,x)\sim(0,T(x))\}
\end{equation*}
defined by the projection
of the horizontals $\intervalleff{0}{1}\times\{x\}$.
Then $(S,\mathcal{F})$ is clearly topologically
conjugated to the suspension of the first-return map
$P^\gamma_{\mathcal{F}}$,
and the dynamical properties of these two dynamical systems are the same.
For instance $\mathcal{F}$ is \emph{minimal}
(namely has all of its leaves
dense) if and only if $P^\gamma_{\mathcal{F}}$ is \emph{minimal}
(namely has all of its orbits
dense).
Moreover if two foliations $\mathcal{F}_1$ and $\mathcal{F}_2$
on $S$ are conjugated by a homeomorphism isotopic to the identity,
then they admit freely homotopic sections on which
their first-return maps
will be topologically conjugated.
We recall that two circle homeomorphisms $f$ and $g$
are \emph{topologically conjugated}
if there exists a homeomorphism $\varphi$ of $\Sn{1}$
such that $g=\varphi\circ f\circ\varphi^{-1}$.
\par In conclusion, any topological conjugacy invariant
of circle homeomorphisms
will yield an isotopy invariant for foliations.
We now define the only such invariant.

\subsection{Circle homeomorphisms and rotation numbers}
We denote by $x\in\R\mapsto [x]\in\Sn{1}=\R/\Z$
the canonical projection
onto the circle, and
by $R_\theta\colon x\in\Sn{1}\mapsto x+\theta\in\Sn{1}$
the \emph{rotation by $\theta\in\Sn{1}$}.
We refer for instance to
\cite[\S 1.1 \& 2.1]{de_faria_dynamics_2022} and
\cite[\MakeUppercase{\romannumeral 1}.1]{de_melo_one-dimensional_1993}
for a proof of the following classical results
due to Henri Poincaré \cite{poincare_sur_1885}.
\begin{propositiondefinition}[Poincaré]\label{propositiondefinition-nombrerotation}
Let $f\in\Homeo^+(\Sn{1})$ be an orientation-preserving homeomorphism
of the circle.
\begin{enumerate}
 \item For any lift $F$ of $f$, the limit
  $\tau(F)\coloneqq\underset{n\to\pm\infty}{\lim}\frac{F^n(x)-x}{n}$
 exists for any $x\in\R$ and is independent of $x$.
 If $G=F+d$ is another lift of $f$ ($d\in\Z$),
 then $\tau(G)=\tau(F)+d$, and
 \begin{equation*}\label{equation-definitionnombrerotation}
\rho(f)\coloneqq[\tau(F)]\in\Sn{1}
 \end{equation*}
is thus a well-defined point called the \emph{rotation number} of $f$.
\item $\rho(f)$ is invariant under topological conjugacies.
The rotation number of any orientation-preserving homeomorphism
$g$ of an oriented topological circle is thus well-defined
by the relation $\rho(g)\coloneqq\rho(g_0)\in\Sn{1}$,
with $g_0\in\Homeo^+(\Sn{1})$ conjugated to $g$
by an orientation-preserving map.
\item If $f\in\Homeo^+(\Sn{1})$ is minimal,
then it is topologically conjugated to the rotation $R_{\rho(f)}$.
\end{enumerate}
\end{propositiondefinition}

Recall that any half-line $l\in\mathbf{P}^+(\R^2)$
induces an oriented \emph{linear foliation}
$\mathcal{F}_l$ on $\Tn{2}$
defined by $\mathcal{F}_l[x]=[x+l]$ for any $[x]\in\Tn{2}$
(with $x\in\R^2\mapsto [x]\in\Tn{2}=\R^2/\Z^2$ the canonical projection).
According to Proposition \ref{propositiondefinition-nombrerotation}.(3),
if a foliation $\mathcal{F}$ of $\Tn{2}$ admits a section $\gamma$
on which the first return map has rotation
number $[\theta]\in\Sn{1}$,
then $\mathcal{F}$ is topologically conjugated
to the linear foliation $\mathcal{F}_{\R(1,\theta)}$.
However the description of $\mathcal{F}$ up to \emph{isotopy}
demands more than the rotation number of its
first-return maps.
\par Let indeed $\gamma'$ be a simple closed curve
freely homotopic to $\gamma$
and disjoint from it,
and $D$ be a positive Dehn twist
around $\gamma'$ whose support is
disjoint from $\gamma$.
Then the first-return map of
$D_*\mathcal{F}$ on $\gamma$ is equal to the one of
$\mathcal{F}$,
although $D_*\mathcal{F}$ is \emph{not} isotopic to
$\mathcal{F}$.
This motivates the introduction of a finer invariant,
a ``global version of the rotation number'' which will
detect the action
of the mapping class group of $\Tn{2}$.

\subsection{Asymptotic cycles}
Originally introduced by
Schwartzman in \cite{schwartzman_asymptotic_1957}
for topological flows of closed manifolds $M$,
the notion of \emph{asymptotic cycle} fulfills this role.
It associates to any suitable orbit $O$ of the flow
the ``best approximation
of $O$ by a closed loop in homology''.
This notion has a natural projective counterpart
for an oriented topological foliation $\mathcal{F}$
of $\Tn{2}$
that we now quickly describe,
referring to \cite{schwartzman_asymptotic_1957,yano_asymptotic_1985}
for more details.
We consider an auxiliary smooth Riemannian metric $g$ on $\Tn{2}$
and its induced distance
$d_{\mathcal{F}}$ on the leaves of $\mathcal{F}$.
For $x\in \Tn{2}$ and $T\in\R$
we denote by $\gamma_{T,x}$ the closed curve of $\Tn{2}$
obtained by
first following
$\mathcal{F}(x)$ from $x$ to the unique point $y\in\mathcal{F}(x)$
such that $d_{\mathcal{F}}(x,y)=T$, and then closing the curve
by following the minimal geodesic of $g$
from $y$ to $x$.
Following \cite{schwartzman_asymptotic_1957,yano_asymptotic_1985},
the \emph{projective asymptotic cycle} of $\mathcal{F}$ at $x$ is
then defined as the half-line
\begin{equation}\label{equation-asymptoticcycle}
 A_{\mathcal{F}}(x)\coloneqq\R^+\left(\underset{T\to+\infty}{\lim}\frac{1}{T}[\gamma_{T,p}]\right)\in
 \mathbf{P}^+(\Homologie{1}{\Tn{2}}{\R})
\end{equation}
in the first homology group of $\Tn{2}$,
if this limit exists and does not vanish.
This cycle is by definition
constant on leaves,
does not depend on the auxiliary Riemannian metric,
and is moreover natural with respect to any homeomorphism $f$:
\begin{equation*}\label{equation-cycleasymptotiquenaturel}
 A_{f_*\mathcal{F}}(f(x))=f_*(A_{\mathcal{F}}(x)).
\end{equation*}
In particular, any homeomorphism isotopic to the identity acts trivially
on projective asymptotic cycles
(see \cite[Theorem p.275]{schwartzman_asymptotic_1957}).
While the asymptotic cycles have \emph{a priori}
no reason to exist at any point,
they are easily described for foliations of the torus
by the following result
which is a reformulation of
\cite[Theorem 6.1 and Theorem 6.2]{yano_asymptotic_1985}.
We identify henceforth $\Homologie{1}{\Tn{2}}{\R}$ with $\R^2$
through the isomorphism induced by the canonical
projection $\R^2\to\Tn{2}$.
\begin{proposition}[\cite{yano_asymptotic_1985}]
\label{proposition-cyclesasymptotiques}
Let $\mathcal{F}$ be an oriented topological
foliation of $\Tn{2}$
which is the suspension of a circle homeomorphism.
\begin{enumerate}
 \item $A_{\mathcal{F}}(x)$ exists at any $x\in\Tn{2}$.
 It is moreover
 constant on $\Tn{2}$ and will be denoted by $A(\mathcal{F})$.
 \item For any $l\in\mathbf{P}^+(\R^2)$,
 $A(\mathcal{F}_l)=l$.
\end{enumerate}
\end{proposition}

Asymptotic cycles play their expected role
of ``global version of the rotation number'',
precisely formulated by the following result
which is folklore in the literature.
\begin{lemma}\label{lemma-cyclesrotation}
 Let $\mathcal{F}_1$ and $\mathcal{F}_2$ be two
 oriented topological foliations of $\Tn{2}$.
 Then $A(\mathcal{F}_1)=A(\mathcal{F}_2)$
 if, and only if
 for any respective sections $\gamma_1$ and $\gamma_2$
 of $\mathcal{F}_1$ and $\mathcal{F}_2$
 which are freely homotopic, we have
 $\rho(P^{\gamma_1}_{\mathcal{F}_1})=\rho(P^{\gamma_2}_{\mathcal{F}_2})$.
\end{lemma}

Using the previous Lemma,
one easily obtains the following classification
result usually attributed to Poincaré.
\begin{proposition}\label{proposition-classificationfeuilletagesminimaux}
Let $\mathcal{F}$ be a minimal oriented topological
foliation of $\Tn{2}$.
Then $\mathcal{F}$ is conjugated
to the linear foliation defined by its asymptotic cycle
$A(\mathcal{F})\in\mathbf{P}^+(\R^2)$,
by a homeomorphism isotopic to the identity.
\end{proposition}

\section{Proof of Theorem \ref{theoremintro-bifoliations}}
\begin{definition}\label{definition-transverse}
 A pair $(\Falpha,\Fbeta)$ of topological foliations
of $\Tn{2}$ is said \emph{transverse}
if for any $p\in\Tn{2}$ there exists
a connected open neighborhood $U$ of $p$,
two open intervals $I,J\subset\R$ containing $0$
and a homeomorphism
$\varphi\colon U\to I\times J$,
which sends every plaque of
$\Falpha$ (respectively $\Fbeta$) in $U$
to a horizontal interval $I\times\{y\}$
(resp. vertical interval $\{x\}\times J$).
Such a homeomorphism will be called a
\emph{simultaneous foliated chart} of $(\Falpha,\Fbeta)$.
\end{definition}

We now prove Theorem \ref{theoremintro-bifoliations}.
According to Proposition \ref{proposition-classificationfeuilletagesminimaux},
we can assume without lost of generality that
$\Fbeta$ is a linear foliation.
We denote by
$\Falphazero$ the linear foliation defined
by $A(\Falpha)\in\mathbf{P}^+(\Homologie{1}{\Tn{2}}{\R})$,
fix $p\in\Tn{2}$ and denote
$\Fa\coloneqq\Falpha(p)$ and
$\Fz\coloneqq\Falpha^0(p)$.

\par \textbf{(a) Flowing along $\Fbeta$ on $\Falpha(p)$.}
Let $U$ be the domain of a simultaneous
foliated chart of $(\Falpha,\Fbeta)$ around $p$.
Then for any point $x$ in the plaque $P_{\Falpha}^U(p)$ of $p$,
there exists a unique point $\phi(x)\coloneqq
P_{\Falphazero}^U(p)\cap P^U_{\Fbeta}(x)$
belonging both to the plaque of $p$ for $\Falphazero$
(defined as the connected component of
$\mathcal{F}_\alpha^0(p)\cap U$
containing $p$)
and to the plaque of $x$ for $\Fbeta$ in $U$.
Let us now assume
$\phi$ to be defined on an interval $I=\intervalleff{a}{b}$
of $\Fa$ with values in $F_\alpha^0$.
Let $\gamma\colon\intervalleff{0}{1}\to\Fbeta(b)$
denote a continuous parametrization of the (unique)
interval of $\Fbeta(b)\setminus\{b,\phi(b)\}$
whose closure contains both $b=\gamma(0)$ and $\phi(b)=\gamma(1)$.
Then with $U$ and $V$ the domains of respective
foliated charts of $\Falpha$ at $b$ and of $\Falpha^0$ at $\phi(b)$,
we extend $\phi$ on $P^U_{\Falpha}(b)$ to be equal to
the holonomy map
of $\Fbeta$ along $\gamma$ from
$P^U_{\Falpha}(b)$ to $P^V_{\Falpha^0}(\phi(b))$
(which is well-defined, possibly diminishing $U$).
The subset of $\Fa$ on which $\phi$ is uniquely defined in this way
is thus non-empty, open and closed, hence equal to $\Fa$.

\par \textbf{(b) Extending $\phi$ to $\Tn{2}$.}
The only possible candidate for a continuous extension
(denoted in the same way) of $\phi$ is of course
\begin{equation}\label{equation-extensionvarphi}
 \phi(\lim x_n)\coloneqq\lim \phi(x_n)
\end{equation}
for any converging sequence $x_n\in F_\alpha$.
Our main task is thus to show first that
for any converging sequence $x_n\in F_\alpha$,
the sequence $\phi(x_n)$ converges,
and to show moreover that for two converging sequences
$x_n^1$ and $x_n^2$ in $F_\alpha$ having the same limit  $x$,
we have $\lim \phi(x_n^1)=\lim \phi(x_n^2)$.
Let $y_n^{i}$ ($i=1$ or $2$) denote the first intersection point
of $\Falpha(x_n^{i})$ with $\Fbeta(x)$,
and let us temporarily admit
that $\phi(y_n^{1/2})$ are convergent and
have the same limit.
For any large enough $n$,
let $J_{+,n}^i$ denote the segment of leave of $\Fbeta(x)$
from $y_n^i$ to $\phi(y_n^i)$,
$I_{-,n}^i$ denote the segment of leave of $\Falpha(x_n^i)$
from $x_n^i$ to $y_n^i$,
and $I_{+,n}^i$ denote the segment of leave of $\Falpha(\phi(x_n^i))$
from $\phi(x_n^i)$ to $\phi(y_n^i)$.
Then by definition of $\phi$,
there exists a unique segment $J_{-,n}^i$ of $\Fbeta(x_n^i)$
such that there exists an embedded rectangle $\mathcal{R}$ of boundary
$I_{-,n}^i\cup J_{+,n}^i\cup I_{+,n}^i\cup J_{-,n}^i$,
and foliated by segments of $\Falpha$ and $\Fbeta$.
This shows that $\phi(x_n^i)=J_{-,n}^i\cap I_{+,n}^i$
is entirely described by $x_n^i$ and $\phi(y_n^i)$,
that the convergence of $x_n^{1/2}$ and $\phi(y_n^{1/2})$ together
imply the one of $\phi(x_n^{1/2})$,
and that the equalities
$\lim x_n^1=\lim x_n^2$ and
$\lim \phi(y_n^1)=\lim \phi(y_n^2)$
imply $\lim \phi(x_n^1)=\lim \phi(x_n^2)$.
\par To prove that
$\lim \phi(y_n^1)$ and $\lim \phi(y_n^2)$ exist and are equal,
we can assume without lost of generality that $y_n^{1/2}$
converges to $x$ from above on $\Fbeta(x)$, and
that the sequence $y_n^{1/2}$ is also increasing
on the oriented topological line $\Falpha(p)$.
Let $a$ denote a (long) segment
of $\Fbeta(x)$ containing $x$, $y_n^1$, $y_n^2$, $\phi(y_n^1)$
and $\phi(y_n^2)$
for any large enough $n$.
We can then close $a$
with a (short) path $b$
to obtain a simple closed curve $\gamma=a\cdot b$ transverse to $\Falpha$.
Furthermore since $\Falpha$ and $\Falpha^0$ are isotopic
and both transverse to $\Fbeta$,
there exists a (short) path $b^0$ and
a segment $a^0$ of $\Fbeta(x)$
with $a\cap a^0$ connected and
containing $x$, $y_n^1$, $y_n^2$, $\phi(y_n^1)$ and $\phi(y_n^2)$
for any large enough $n$,
such that $\gamma^0=a^0\cdot b^0$ is a simple closed curve transverse
to $\Falpha^0$ and homotopic to $\gamma$.
Then with $P$ and $P_0$ the respective first-return maps of
$\Falpha$ and $\Falpha^0$ on $\gamma$ and $\gamma^0$,
we have
\begin{equation}\label{equation-egalitenombrerotationfirstreturn}
 [\theta]\coloneqq\rho(P)=\rho(P_0)\in\Sn{1}=\R/\Z,
\end{equation}
since $A(\Falpha)=A(\Falpha^0)$
and $\gamma$ is homotopic to $\gamma^0$.
Let $y$ be the first intersection point of the oriented leaf
$\Falpha(p)$ with $\gamma$.
Since $y_n^{1/2}$ is increasing
on $\Falpha(p)$,
there exists
two increasing sequences $k_n^1,k_n^2\in\N$ such that
$y_n^{1/2}=P^{k_n^{1/2}}(y)$.
With $y^0$ the first intersection point of
$\Falpha^0(p)$ with $\gamma$,
we have then
$\phi(y_n^{1/2})=P_0^{k_n^{1/2}}(y^0)$
by the very definition of $\phi$.
\par We now make a naive but useful general remark.
Let $f$ be a minimal orientation-preserving circle
homeomorphism
of rotation number $[\theta]\in\Sn{1}$.
Then since $f$ is topologically conjugated to the rotation $R_\theta$,
for any $x\in\Sn{1}$ and $k_n\in\N$,
the sequence $(f^{k_n}(x))$ is converging if, and only if
$[k_n\theta]$ is converging in $\Sn{1}$.
Moreover for $k_n^1,k_n^2\in\N$
such that $(f^{k_n^1}(x))_n$ and $(f^{k_n^2}(x))_n$ are converging,
$\lim f^{k_n^1}(x)=\lim f^{k_n^2}(x)$ if, and only if
$\lim [(k_n^1-k_n^2)\theta]=[0]\in\Sn{1}$.
Since $P$ and $P^0$ have the same rotation numbers
according to
\eqref{equation-egalitenombrerotationfirstreturn},
the convergence of
$y_n^1=P^{k_n^1}(y)$ and $y_n^2=\lim P^{k_n^2}(y)$
and the equality of their limits
is thus equivalent to the convergence of
$\phi(y_n^1)=P^{k_n^1}(y)$ and $\phi(y_n^2)=P^{k_n^2}(y)$
and to the equality of their limits.
Therefore $\lim \phi(y_n^1)=\lim \phi(y_n^2)$,
and we can thus extend $\phi$ as desired, to an application
well-defined on $\Tn{2}$
by the relation \eqref{equation-extensionvarphi}.
\par Note that we incidentally proved that $\phi$
is ``continuous along $F_\alpha$'', in the sense
that the equality \eqref{equation-extensionvarphi}
holds for any sequence $x_n\in F_\alpha$ converging
to a point of $F_\alpha$.
Moreover, we also proved that
$\lim\phi(x_n^1)=\lim\phi(x_n^2)$ implies that $x^1=x^2$
for any two sequences $x_n^1$ and $x_n^2$ in $F_\alpha$
respectively converging to $x^1$ and $x^2$,
hence that $\phi$ is injective by our definition of $\phi$.

\par \textbf{(c) $\phi$ is a homeomorphism.}
Let us first show that $\phi$ is continuous.
For any sequence $x_n\in\Tn{2}$ converging to $x$,
let $x_n^k\in F_\alpha$ denote for any $n$ a sequence
converging to $x_n=\underset{k\to+\infty}{\lim}x_n^k$.
Then by the definition \eqref{equation-extensionvarphi}
of $\phi$, we have
$\phi(x_n)=\underset{k}{\lim}~\phi(x_n^k)$.
Possibly extracting a subsequence of $(x_n^k)_k$,
we can furthermore assume that
\begin{equation}\label{equation-controlexnk}
 d(\phi(x_n^k),\phi(x_n))\leq\frac{1}{k+1}
\end{equation}
for any $n$ and $k$,
with $d$ a distance defining the topology of $\Tn{2}$.
There exists now an increasing sequence $k_n\in\N$
such that
$x_n^{k_n}\in F_\alpha$ converges to $x$,
and we have thus $\phi(x)=\lim\phi(x_n^{k_n})$
by the definition \eqref{equation-extensionvarphi}
of $\phi$.
But since $d(\phi(x_n^{k_n}),\phi(x_n))\leq\frac{1}{k_n+1}$
according to \eqref{equation-controlexnk}, we obtain
$\phi(x)=\lim~\phi(x_n^{k_n})=\lim~\phi(x_n)$
which proves the continuity of $\phi$.
\par Since $\phi$ is an injective and continuous map
defined from the topological surface $\Tn{2}$ to itself,
the Invariance of domain theorem of Brouwer
shows that $\phi(\Tn{2})$ is an open (and non-empy) subset
of $\Tn{2}$.
Since it is also closed by compactness of $\Tn{2}$,
we eventually have $\phi(\Tn{2})=\Tn{2}$
by connectedness of $\Tn{2}$,
hence $\phi$ is a homeomorphism of $\Tn{2}$ as desired.
\par \textbf{(d) $\phi$ preserves $\Fbeta$.}
By construction, the restriction of $\phi$ to $F_\alpha$
consists by flowing along leaves of $\Fbeta$.
We have therefore $\phi(y)\in\Fbeta(\phi(x))$
for any $x\in F_\alpha$ and $y\in F_\alpha\cap\Fbeta(x)$.
By continuity of the foliations and of $\phi$, this relation
extends to $\Tn{2}$, hence $\phi$ preserves $\Fbeta$.

\par \textbf{(e) $\phi$ is isotopic to $\id$ rel $p$.}
Moreover $\phi$ fixes $p$ by construction,
and acts trivially on $\pi_1(\Tn{2},p)$.
Indeed, let $(a,b)$ be a pair of simple closed curves based at $p$
which are the concatenations $a=a_\alpha a_\beta$ and
$b=b_\alpha b_\beta$ of segments $a_{\alpha/\beta}$
(respectively $b_{\alpha/\beta}$) of leaves of
$\mathcal{F}_{\alpha/\beta}$,
and whose homotopy classes define a basis $([a],[b])$
of $\pi_1(\Tn{2},p)$.
Since $\phi$ is defined on $\Falpha(p)$ by flowing along leaves of
$\Fbeta$, the paths $a_\alpha$ (respectively $b_{\alpha}$)
are isotopic to
$\phi\circ a_{\alpha}$ (resp. $\phi\circ b_{\alpha}$),
their arrival point flowing along $\Fbeta(p)$.
But the paths $a_{\beta}$ (resp. $b_{\beta}$)
are for the same reason isotopic to
$\phi\circ a_{\beta}$ (resp. $\phi\circ b_{\beta}$),
their departure point flowing along $\Fbeta(p)$,
and therefore $a$ (resp. $b$) is isotopic to
$\phi\circ a$ (resp. $\phi\circ b$).
Since $\phi$ acts trivially on $\pi_1(\Tn{2},p)$,
it is isotopic to the identity
by a classical result of Epstein in
\cite{epstein_curves_1966}
(see also \cite[Proposition 1.6 and Theorem 2]{beguin_fixed_2020}).

\par \textbf{In conclusion},
$\phi$ is a homeomorphism isotopic to $\id$ relatively to $p$,
preserving $\Fbeta$ and redressing $\Falpha$ on $\Falpha^0$,
\emph{i.e.} is a conjugation of $(\Falpha,\Fbeta)$
with $(\Falpha^0,\Fbeta)$.
This concludes
the proof of Theorem \ref{theoremintro-bifoliations}.

\section{Proof of Corollary \ref{corollaryintro-automorphismesbifoliations}}
We conclude this note with a proof of
Corollary \ref{corollaryintro-automorphismesbifoliations}.
According to Theorem \ref{theoremintro-bifoliations},
the minimal bi-foliation preserved by $f$ is conjugated to
a linear bi-foliation
by a homeomorphism isotopic to the identity.
We can thus assume without
lost of generality that $f$ preserves the bi-foliation
defined by two transverse irrational lines $l$ and $l'$
of the plane, and
is isotopic to the identity relatively to $[0]\in\Tn{2}$.
Let
\begin{equation*}\label{eqation-F}
 F(x,y)=(F_1(x,y),F_2(x,y))
\end{equation*}
denote the lift of $f$ to $\R^2$ which fixes the origin.
According to
\cite[Lemme 4]{ghysStabiliteConjugaisonDifferentiable1980},
we have then
\begin{equation*}\label{equation-definitionxprimeyprimeg1et2}
\begin{cases}
 F_2(x,y)=\delta F_1(x,y)+a(y-\delta x)+b \\
 F_2(x,y)=\delta' F_1(x,y)+a'(y-\delta' x)+b'
 \end{cases}
\end{equation*}
with $\delta$ and $\delta'$ the respective slopes of $l$ and $l'$.
A direct computation yields then
\begin{equation*}\label{equation-definitionxprimeyprimeg1et2}
\begin{cases}
 F_1(x,y)=\frac{1}{\delta-\delta'}\left(
 (a\delta-a'\delta')x+(a'-a)y+(b'-b)
 \right) \\
 F_2(x,y)=\frac{1}{\delta-\delta'}\left(
 \delta\delta'(a-a')x+(\delta a'-\delta'a)y+(\delta b'-\delta'b)
 \right).
 \end{cases}
\end{equation*}
In other words:
$F(x,y)=M(x,y)+\frac{1}{\delta-\delta'}(b'-b,\delta b'-\delta'b)$
with
\begin{equation*}\label{equation-M}
 M=\frac{1}{\delta-\delta'}
 \begin{pmatrix}
    a\delta-a'\delta' & a'-a \\
    \delta\delta'(a-a') & \delta a'-\delta'a
   \end{pmatrix}.
\end{equation*}
Since $F(0,0)=(0,0)$ and $\delta\neq\delta'$, we have $b=b'=0$.
Moreover $M\in\mathrm{GL}_2(\Z)$ since $F$
induces the map $f$ on $\Tn{2}$,
and the action of $f$ on $\piun{\Z^2}$ is thus given by the matrix $M$.
Since $f$ is isotopic to the identity, this implies $M=\id$
and therefore $F=\id$, which shows our claim and concludes the proof
of Corollary \ref{corollaryintro-automorphismesbifoliations}.

\bibliographystyle{alpha}
\bibliography{bifoliationsT2-biblio}

\end{document}